\documentclass[11pt]{amsart}
\usepackage{amsfonts,amssymb,amsmath}
\usepackage[all]{xy}
\begin{document}

\newcommand{\Ext}{{\rm Ext}}
\newcommand{\F}{{\mathbb F}}
\newcommand{\Q}{{\mathbb Q}} 
\newcommand{\T}{{\mathbb T}}
\newcommand{\Z}{{\mathbb Z}}
\newcommand{\hZ}{{\widehat{\mathbb Z}}}
\newcommand{\Zo}{{{\mathbb Z}_0}}
\newcommand{\R}{{\mathbb R}}
\newcommand{\cE}{{\mathcal E}}
\newcommand{\bE}{{\mathbb E}}
\newcommand{\tE}{{\widetilde{\mathcal E}}}
\newcommand{\oD}{{\overline{\Delta}}}
\newcommand{\A}{{\mathbb A}}
\newcommand{\Hom}{{\rm Hom}}
\newcommand{\oh}{{\mathfrak o}}

\title{ A topological group of extensions of $\Q$ by $\Z$}
\author{Jack Morava}
\address{The Johns Hopkins University,
Baltimore, Maryland 21218}
\email{jack@math.jhu.edu}

\subjclass{11R56, 14G40, 18G15}
\date{13 Oct 2013}
\begin{abstract}{The group of extensions (as in the title), endowed with 
something like a connection at Archimedean infinity, is isomorphic to the 
ad\'ele-class group of $\Q$: which is a topological group with interesting 
Haar measure.}
\end{abstract}\bigskip

\maketitle \bigskip

{\it For Mike Boardman and Takashi Ono: dear friends and colleagues}\bigskip

\noindent {\bf 1.1} An extension of one abelian group $A$ (ie, a $\Z$-module), by another ($B$) is an
exact sequence
\[
\xymatrix{
\cE : 0 \ar[r] & B \ar[r]^i & E \ar[r]^j & A \ar[r] & 0 \;;}
\]
the exact functor
\[
C \mapsto C\otimes_\Z \R := C_\R
\]
associates to $\cE$ an extension
\[
\xymatrix{
\cE_\R : 0 \ar[r] & B_\R \ar[r]^{i_\R} & E_\R \ar[r]^{j_\R} & A_\R \ar[r] & 0 }
\]
of real vector spaces, which necessarily splits. This note is concerned with extensions $\cE$ as 
above, which have been {\bf rigidified} by the choice 
\[
s_\cE : A_\R \to E_\R
\]
of a splitting of $\cE_\R$, ie a homomorphism such that $j_\R \circ s_\cE = 1_{A\R}$. I will refer to 
the pair $\tE := (\cE,s_\cE)$ as an extension of $\Z_0$-modules.\bigskip

\noindent
A {\bf congruence} $\alpha$ of two extensions $\cE, \; \cE'$ of $A$ by $B$ is a commutative
diagram
\[
\xymatrix{
\ 0 \ar[r] & B \ar[d]^{1_B} \ar[r]^i & E \ar[d]^\alpha \ar[r]^j & A \ar[d]^{1_A} \ar[r] & 0 \\
0 \ar[r] & B' \ar[r]^{i'} & E' \ar[r]^{j'} & A \ar[r] & 0 \;,}
\]
cf eg [6 III \S1]. Let us say that a congruence $\alpha: \tE \equiv \tE'$ of rigidified extensions 
is a congruence $\alpha: \cE \equiv \cE'$ of their underlying extensions of $\Z$-modules, which is 
compatible with their splittings in the sense that the diagram 
\[
\xymatrix{E_\R \ar[d]^{\alpha_\R}  & \ar[l]_{s_\cE} A_\R \ar[d]^{1_{A\R}} \\
E'_\R & \ar[l]_{s_{\cE'}} A_\R }
\]
commutes.\bigskip

\noindent
Congruence classes of such rigidified extensions define an abelian group-valued bifunctor $\Ext_\Zo(A,B)$, 
with a straightforward generalization of Baer sum (as we shall check below) as composition. Forgetting 
the splitting data defines an epimorphism
\[
\Ext_\Zo(A,B) \to \Ext_\Z(A,B) \;.
\]
\bigskip

\noindent
{\bf 1.2.1 Proposition} This forgetful homomorphism fits in an exact sequence
\[
0 \to \Hom_\R(A_\R,B_\R) \to \Ext_{\Z_0}(A,B) \to \Ext_\Z(A,B) \to 0 \;;
\]
in particular, the exact sequence
\[
0 \to \R \to \Ext_{\Z_0}(\Q,\Z) \to \Ext_\Z(\Q,\Z) \to 0
\]
is isomorphic to the sequence
\[
0 \to \R := \Sigma_0 \to \Sigma \to (\hZ \otimes \Q)/\Q \to 0
\]
defined by the inclusion of the path-component $\Sigma_0$ of the identity in the solenoid
\[
o \to \Q \to \A_\Q := \R \times (\hZ \otimes \Q) \to \Sigma \cong \Hom(\Q,\T) \to 1
\]
(ie the Pontrjagin dual of the rational numbers: which is connected but not path-connected).\bigskip

\noindent
{\bf 1.2.2} More generally, if $\oh_K$ is the ring of algebraic integers in a number field $K$, then
\[
\Ext_{\Z_0}(\Q,\oh_K) \; \cong \; \A_K/K
\]
is naturally isomorphic to the ad\'ele-class group of $K$ [3 \S 14, 5 \S 5.3]: a compact topological group 
(a product of $[K:\Q]$ copies of $\Sigma$) with canonical Haar measure whose square equals the absolute 
value of the discriminant $D_{K/\Q}$ of $K$ over $\Q$. 
\bigskip

\noindent
{\bf Remark} The injective resolution 
\[
0 \to \Z \to \Q \to \Q/\Z \to 0
\]
identifies
\[
\Ext_\Z(\Q,\Z) \cong (\hZ \otimes \Q)/\Q
\]
with the cokernel of the map defined by tensoring the (dense) inclusion
\[
\Z \to \hZ = \prod \Z_p
\]
(of the integers into the product of the profinite integers) with $\Q$: which is superfluous, since
$\hZ/\Z$ is uniquely divisible and is thus already a $\Q$-vector space, whose natural [?] topology
is then {\bf indiscrete}. See [2 Theorem 25] for an account (which motivated this note) of this 
classical group of extensions. \bigskip

\noindent
{\bf 1.3.1 Proof:} We need first to define the Baer sum of two rigidified extensions $\tE_0, \; \tE_1$ 
as the quotient $\tE_+$ of the pullback $\bE_+$
\[
\xymatrix{
B \ar[dr] \ar[drr] \ar[ddr] & {} & {} \\
{} & \bE_+ \ar[d] \ar[r] & E_1 \ar[d]^{j_1} \\
{} & E_0 \ar[r]^{j_0} & A }
\]
by the image of the map
\[
b \mapsto i_\Delta(b) = (i_0(b),-i_1(b)) : B \to \bE_+ \;.
\]
The resulting sum is split, after tensoring with $\R$, by
\[
A \ni a \mapsto (s_0(a),s_1(a)) \in E_{+ \R} = \bE/({\rm image}
 \; i_\Delta) \;.
\] 

\noindent
To check the first assertion of the proposition, note that if two extensions $\tE,\; \tE'$ of $\Z_0$-modules
have congruent underlying extensions of $\Z$-modules, then tensoring those extensions with $\R$ defines an 
isomorphism
\[
\xymatrix{
0 \ar[r] & B_\R \ar[d] \ar[r]^i & E_\R \ar[d]^{\alpha_\R} \ar[r]^j &
A_\R \ar[d] \ar[r] & 0 \\
0 \ar[r] & B_\R \ar[r]^{i'} & E'_\R \ar[r]^{j'} & A_\R \ar[r]  & 0 }
\]
of extensions of real vector spaces, with splittings $s,s'$.
If
\[
\rho' := s' - s \circ \alpha_\R : A_\R \to E'_\R
\]
then $j' \circ \rho' = 0$, so the image of $\rho'$ lies in the image of $i'$, defining
\[
(i')^{-1} \circ \rho' \in \Hom_\R(A_\R,B_\R) \;.
\]
On the other hand such homomorphisms act freely on the $\Z_0$-module extensions of $A$ by
$B$: $\Hom_\R(A_\R,B_\R)$ is the group, under Baer sum, defined by splittings of the exact
sequence
\[
0 \to B_\R \to B_\R \oplus A_\R \to A_\R \to 0 \;.
\]

\noindent
To prove the second assertion of the proposition we construct a homomorphism
\[
\oD : \Ext_{\Z_0}(\Q,\Z) \to \Hom(\Q,\T)
\]
as follows: given a diagram
\[
\xymatrix{
0 \ar[r] & \Z \ar[d] \ar[r]^i & E \ar[d] \ar[r]^j & \Q \ar[d] \ar[r] & 0 \\
0 \ar[r] & \R \ar[r]^{i_\R} & E_\R \ar[r]_{j_\R} & \ar[l]_s \R \ar[r] & 0 \;,}
\]
let
\[
\Delta_0(e) := e_\R - s(j(e)_\R) : E \to E_\R
\]
(where $x_\R := x \otimes 1_\R$); then $j_\R \circ \Delta_0 = 0$, so the image of $\Delta_0$ lies
in the image of $i_\R$, defining
\[
\Delta := i_\R^{-1} \circ \Delta_0 : E \to \R \;.
\]
But now $\Delta \circ i : \Z \to E \to \R$ is the usual inclusion, so $\Delta$ induces a
homomorphism
\[
[\oD : \Q = E/\Z \to \R/\Z] \in \Hom(\Q,\T) \;\dots
\]
\bigskip

\noindent
{\bf 1.3.2} Here is an example:
\[
\xymatrix{
0 \ar[r] & \Z \ar[d] \ar[r]^<<<<i & {\Z_{(p)} \oplus Z[p^{-1}]} \ar[d] \ar[r]^<<<<j & \Q \ar[d] \ar[r] & 0 \\
0 \ar[r] & \R \ar[r] & \R \oplus \R \ar[r] & \ar[l]_\sigma \R \ar[r] & 0 }
\]
with $i(k) = (k,k), \; j(u,v) = u - v$, and $\sigma(x) = ((s+1)x,sx)$ for some $s \in \R$.\bigskip

\noindent
[{\bf Check} that $j$ is onto, ie that $q \in \Q$ equals $u - v$ for some $u \in                                    
\Z_{(p)}$ and $v \in \Z[p^{-1}$]: \bigskip

\noindent
Let $q = q_0p^{-n}/q_1$ with $p \nmid q_0,q_1$ and $n > 0$ (otherwise the claim is
immediate, with $v = 0$). Then $(p^n,q_1) = 1$ implies the existence of $\alpha, \beta$ such
that
\[
\alpha p^n + \beta q_1 = 1 \;,
\]
so if $a = q_0 \alpha, \; b = - q_0 \beta$ then $ap^n - bq_1 = q_0$ and hence
\[
q = \frac{ap^n - bq_1}{q_1 p^n} = \frac{a}{q_1} - \frac{b}{p^n} \;.
\]

\noindent
Now we have $\Delta_0(u,v) = (u,v) - \sigma (u-v) = i(v - s(u-v))$, so
\[
\Delta (u,v) = (s+1)v - su \; {\rm mod} \; \Z \;,
\]
eg $\oD(q) = - sq + q_0\beta p^{-n} \in \T$. \bigskip

\newpage

\noindent
{\bf 1.3.3} More generally, $\Ext_{\Z_0}(\Q,\Z^k)$ is isomorphic to the product $\Sigma^k \cong 
(\A_\Q/\Q)^k$. The ring $\oh_K$ of integers in a number field, however, has more structure,
which endows its ad\'ele-class group $\A_K/K \cong (\A_\Q/\Q)^k$ with a Haar measure normalized 
[8 V \S 4 Prop 7] as asserted in the proposition.\bigskip  

\noindent
{\bf 1.4 Remarks} I suppose the proposition above has a natural reformulation in
Arakelov geometry [4 \S 1], ie in terms of $\Z$-modules $A,B$ endowed with positive-definite
inner products on $A_\R,B_\R$; but I don't know anything about Arakelov geometry\begin{footnote}{Indeed 
this note is essentially a footnote to work of Bost and K\"unnemann [9]: more thanks to ChD for 
this reference!}\end{footnote}. \bigskip

\noindent
It also seems plausible that 
\[
\Ext_{\F[t]}(\F(t),\F[t]) \cong \A_{\F(t)}/\F(t)
\]
for finite fields $\F$; but its analog of Haar measure seems to be more like an 
invariant one-form [1].
\bigskip

\noindent
{\bf 1.5 Acknowledgements} I owe A. Salch [6 \S 3] thanks, for bringing these matters to my attention. 
This note grew out of conversations with Ch. Deninger and R. Meyer at the September 2013 Oberwolfach 
workshop on noncommutative geometry, and I would like to thank them both for their interest and help.\bigskip

\bibliographystyle{amsplain}

\begin{thebibliography}{99}


\bibitem[1]{1} J V Armitage, On a theorem of Hecke in number fields and 
function fields, Invent. Math. 2 (1967) 238 - 246

\bibitem[2]{2} J M Boardman, Some common Tor and Ext groups, available at 
\[
{\tt  math.jhu.edu/ \sim jmb/note/torext.pdf}
\]

\bibitem[3]{3} J W S Cassels, in {\bf Algebraic number theory}, ed Cassels \& Fr\"ohlich

\bibitem[4]{4} M Kapranov, O Schiffmann, E Vasserot, The spherical Hall algebra of Spec $\Z$, available
at {\tt arXiv:1202.4073}

\bibitem[5]{5} H Koch, {\bf Algebraic number theory}, Encyclopaedia of Math. Sci. 62, Number Theory II,
Springer (1992)

\bibitem[6]{6} S Mac Lane, {\bf Homology}, Springer Grundlehren 114 (1963)

\bibitem[7]{7} A Salch, Grothendieck duality under Spec $\Z$, available at {\tt arXiv:1012.0110}

\bibitem[8]{8} A Weil, {\bf Basic number theory}, Springer Grundlehren 144 (1967) \bigskip

\bibitem[9]{9} J B Bost, K K\"unnemann, Hermitian vector bundles and extension groups on arithmetic schemes I\dots,
Adv. Math. 223 (2010) 987 - 1106, {\tt arXiv:math/0701343}

\end{thebibliography}

\end{document}